\documentclass[12pt]{amsart}
\usepackage{}
\usepackage{amsmath}
\usepackage{amsfonts}
\usepackage{amssymb}
\usepackage[all]{xy}           
\usepackage{bm}
\usepackage{bbding}
\usepackage{txfonts}
\usepackage{amscd}
\usepackage{xspace}
\usepackage[shortlabels]{enumitem}
\usepackage{ifpdf}

\ifpdf
  \usepackage[colorlinks,final,backref=page,hyperindex]{hyperref}
\else
  \usepackage[colorlinks,final,backref=page,hyperindex,hypertex]{hyperref}
\fi
\usepackage{tikz}
\usepackage[active]{srcltx}

\makeatletter

\newtheorem{thm}{Theorem}[section]

\newtheorem{cor}[thm]{Corollary}
\newtheorem{pro}[thm]{Proposition}

\theoremstyle{definition}
\newtheorem{rmk}[thm]{Remark}
\newtheorem{defi}[thm]{Definition}

\newcommand{\nc}{\newcommand}
\newcommand{\delete}[1]{}

\nc{\mlabel}[1]{\label{#1}}  
\nc{\mcite}[1]{\cite{#1}}  
\nc{\mref}[1]{\ref{#1}}  
\nc{\mbibitem}[1]{\bibitem{#1}} 

\delete{
\nc{\mlabel}[1]{\label{#1}{\hfill \hspace{1cm}{\bf{{\ }\hfill(#1)}}}}
\nc{\mcite}[1]{\cite{#1}{{\bf{{\ }(#1)}}}}  
\nc{\mref}[1]{\ref{#1}{{\bf{{\ }(#1)}}}}  
\nc{\mbibitem}[1]{\bibitem[\bf #1]{#1}} 
}

\newcommand {\emptycomment}[1]{}

\delete{
\setlength{\baselineskip}{1.8\baselineskip}

\newcommand{\emptycomment}[1]{}

}

\nc{\calo}{\mathcal{O}}
\nc{\oop}{$\mathcal{O}$-operator\xspace}
\nc{\oops}{$\mathcal{O}$-operators\xspace}
\nc{\mrho}{{\bm{\varrho}}}
\nc{\bfk}{\mathbf{K}}
\nc{\invlim}{\displaystyle{\lim_{\longleftarrow}}\,}
\nc{\ot}{\otimes}
\nc{\CV}{\mathbf{C}}

\newcommand{\lon }{\,\rightarrow\,}
\newcommand{\be }{\begin{equation}}
\newcommand{\ee }{\end{equation}}

\newcommand{\g}{\mathfrak g}
\newcommand{\h}{\mathfrak h}

\newcommand{\huaC}{{\mathcal{C}}}

\newcommand{\huaO}{{\mathcal{O}}}

\newcommand{\half}{\frac{1}{2}}

\newcommand{\Courant}[1]{\left\llbracket  #1\right\rrbracket }

\newcommand{\Id}{{\rm{Id}}}

\newcommand{\br}[1]{   [ \cdot,    \cdot  ]   }

\newcommand{\Hom}{\mathrm{Hom}}

\newcommand{\gl}{\mathfrak {gl}}

\newcommand{\ad}{\mathrm{ad}}

\newcommand{\dglas}{{\rm dgLas}}
\newcommand{\dgla}{{\rm dgLa}}
\newcommand{\gla}{{\rm gLa}}
\newcommand{\sgla}{{\rm sgLa}}
\nc{\opt}{operator\xspace}
\nc{\Opt}{Operator\xspace}
\nc{\oprn}{\theta}
\nc{\Oprn}{\Theta}

\begin{document}

\title[Deformations and homotopy of $\calo$-operators]{Deformations and homotopy of Rota-Baxter operators and $\huaO$-operators on Lie algebras}

\author{Rong Tang}
\address{Department of Mathematics, Jilin University, Changchun 130012, Jilin, China}
\email{tangrong16@mails.jlu.edu.cn}

\author{Chengming Bai}
\address{Chern Institute of Mathematics and LPMC, Nankai University,
Tianjin 300071, China}
\email{baicm@nankai.edu.cn}

\author{Li Guo}
\address{Department of Mathematics and Computer Science,
         Rutgers University,
         Newark, NJ 07102}
\email{liguo@rutgers.edu}

\author{Yunhe Sheng}
\address{Department of Mathematics, Jilin University, Changchun 130012, Jilin, China}
\email{shengyh@jlu.edu.cn}

\date{\today}

\begin{abstract}
This article gives a brief introduction to some recent work on
deformation and homotopy theories of Rota-Baxter operators and
more generally $\calo$-operators on Lie algebras, by means of the
differential graded Lie algebra approach. It is further shown
that these theories lift the existing connection between
$\calo$-operators and pre-Lie algebras to the levels
of deformations and homotopy.
\end{abstract}

\subjclass[2010]{17B37,81R50,17B56,81R12,16T26,17A30,17B62}

\keywords{deformation, homotopy, Maurer-Cartan equation, $\mathcal{O}$-operator, Rota-Baxter operator, pre-Lie algebra}

\maketitle

\vspace{-.9cm}

\allowdisplaybreaks

\section{Introduction}\mlabel{sec:intr}
This a abridged survey of some of the recent developments on
deformations and homotopy of $\calo$-operators and Rota-Baxter
operators in relation with pre-Lie algebras. Our focus will be
on~\mcite{TBGS-1,TBGS-2} to which we refer the reader for details.

\subsection{Classical theories of deformations and homotopy}
\label{ss:deform} Deformation is an omnipresent notion in
mathematics and physics. Roughly speaking, a deformation of an
object equipped with a certain structure is a perturbation of the
object (by a parameter for instance) which has the same kind of
structure. In physics, the idea of deformation is behind the
perturbative quantum field theory. A regularization procedure in
the renormalization of getting rid of a divergency is often a
deformation in practice. Deformation quantization has been studied
under many contexts in mathematical physics~\cite{K1,K2,Ri,Sc}.

The foundational  work of Kodaira and Spencer~\cite{KS} for deformations of complex analytic structures led to its generalization in algebraic geometry and number theory.
In algebra, deformation theory began with the seminal
work of Gerstenhaber~\cite{Ge0,Ge} for associative
algebras, followed by its extension to Lie algebras by
Nijenhuis and Richardson~\cite{NR,NR2}.

Homotopy is a notion closely related to deformation. Starting with the fundamental procedure in topology describing continuously deforming one function to another, the homotopy of an algebraic structure is obtained when the defining relations of the algebraic structure is relaxed
to hold up to a weaker homotopy condition. Thus homotopy is often regarded as deformation from another point of view. This notion abstracted to category theory gives a mathematical context to the fundamental gauge principle in physics that it is more useful to relate two objects by a gauge transformation rather than a strict equality~\mcite{nLab}. The processes can be iterated and give rise to higher homotopies and higher gauge transformations.

The first homotopy construction in pure algebra is the $A_\infty$-algebra of Stasheff, arising from his work on homotopy
characterization of connected based loop spaces~\cite{St63}. This was followed by $L_\infty$-algebras, both having many applications in physics, especially in supergravity and string theory~\cite{St19}. See also~\cite{Ke} for application of homotopy in statistical physics.
After deformations and homotopy of other algebraic structures such as pre-Lie algebras were developed~(see for example \cite{Bu0,CL}), their general theories for algebraic structures in the context of operads were established~\cite{LV,Mar,Ma}.

However, these general theories mostly apply only to connected operads, namely the operads whose space of unary operations is no bigger than the linear span of the identity automorphism.

\subsection{Rota-Baxter operators and \oops}

It is against this background that the characters of our
studies in~\cite{TBGS-1,TBGS-2} took the stage, motivating us to expand the existing literature of algebraic deformations and homotopy. The characters are the
$\calo$-operators and Rota-Baxter operators, with the latter being a special case of the former.

Since this paper is intended to be a brief survey, we will
restrict the weight of the $\calo$-operators and Rota-Baxter
operators to be zero, and limit the discussion on
$r$-matrices to general remarks.

\begin{defi} \label{defi:O} Let $(\g,[\cdot,\cdot])$ be a Lie algebra.
\begin{enumerate}
\item[\rm(i)]
 A linear operator $P:\g\longrightarrow \g$ is called a {\bf Rota-Baxter operator} (of weight $0$) if
\begin{equation} [P(x),P(y)]=P\big([P(x),y]+ [x,P(y)]\big), \quad \forall x, y \in \g.
\label{eq:rbo}
\end{equation}
\item[\rm(ii)]
Let $\rho:\g\longrightarrow\gl(V)$ be a representation of $\g$ on a vector space $V$. An {\bf \oop} on $\g$ with respect to the representation $(V;\rho)$ is a linear map $T:V\longrightarrow\g$ such that
 \begin{equation}
   [Tu,Tv]=T\big(\rho(Tu)(v)-\rho(Tv)(u)\big),\quad\forall u,v\in V.
 \mlabel{eq:defiO}
 \end{equation}
\end{enumerate}
\label{de:conc}
\end{defi}

Note that when $\rho$ is the adjoint representation of $\g$,
Eq.~\eqref{eq:defiO} reduces to Eq.~(\ref{eq:rbo}), which means that a Rota-Baxter operator is an \oop on $\g$ with respect to the adjoint representation. Furthermore, a skew-symmetric $r$-matrix
corresponds to an $\mathcal O$-operator on $\g$ with respect to
the coadjoint representation \cite{Ku}.

The concept of Rota-Baxter operators on associative algebras was
introduced in 1960 by G. Baxter \cite{Ba} in his study of
fluctuation theory in probability. Recently it has found many
applications, including in Connes-Kreimer's
algebraic approach to the renormalization of perturbative quantum
field theory~\cite{CK,Gub}. In the Lie algebra context, a Rota-Baxter operator
(of weight 0) was introduced independently in the 1980s as the
operator form of the classical Yang-Baxter equation, named after
the well-known physicists C. N. Yang and R.~J.~Baxter~\cite{STS}.

To better understand the classical Yang-Baxter equation and
the related integrable systems, the more general notion of an \oop (later also called a relative Rota-Baxter operator or a generalized Rota-Baxter operator) on a Lie algebra was introduced~\cite{Bor,Ku}.

\subsection{Our approach}
Given the importance of Rota-Baxter operators and the more general $\calo$-operators in mathematics and mathematical physics, it is meaningful to study their deformations and homotopies. However, as noted above, the existing general framework of deformation and homotopy for operads does not apply to Rota-Baxter operators and hence $\calo$-operators: the operad of Rota-Baxter (associative or Lie) algebras has the Rota-Baxter operator as the non-identity unary operation and hence is not connected and thus is not covered by the general theories of deformations and homotopy.

Instead, we take our approach by general principles of
deformation and homotopy theories. As stated in~\cite{St19}, the
following ``metatheorem" provides a philosophy (also see the
letter of Deligne~\cite{De}) for a deformation (or a homotopy)
theory:

\begin{quote}
The deformation theory of any mathematical object, e.g., an associative algebra, a complex manifold, etc., can be described starting from a certain differential graded Lie algebra (dgLa) associated to the mathematical object in question.
\end{quote}
Further, the deformations are given as solutions of a ``Master Equation", now known as the Maurer-Cartan equation on the dgLa.

In our constructions of deformation and homotopy theories of Rota-Baxter operators and $\calo$-operators, we also ensure that the theories are compatible with the recently established relationship of the operators with pre-Lie algebras. To be specific, it has been proved that an $\calo$-operator on a Lie algebra naturally gives rise to a pre-Lie algebra~\cite{Bai}, whose deformation and homotopy theories have been established~\cite{Bu0,CL}. Thus a test for our deformation and homotopy theories of $\calo$-operators is whether they extend the above connection of $\calo$-operators with pre-Lie algebras to the contexts of deformations and homotopy. We show that this is indeed the case.

The outline of the paper is as follows. In Section~\ref{sec:def}, we first recall background on differential graded Lie algebras ($\dglas$) and then apply it to define the deformations of $\calo$-operators. We then establish a natural homomorphism from the $\dgla$ defining the deformations of $\calo$-operators to the one defining the deformation of pre-Lie algebras, yielding the desired relationship between the deformation theories. In Section~\ref{sec:hom}, we take a similar approach, but work in the framework of symmetric $\dglas$, for the homotopy of $\calo$-operators including its connection with the homotopy of pre-Lie algebras.

\section{Deformations of $\huaO$-operators}
\label{sec:def}
We first recall some general notions. Let $V=\oplus_{k\in\mathbb Z}V^k$ be a $\mathbb Z$-graded vector space. We denote by $S(V)$ the symmetric  algebra  of $V$, i.e.
$
S(V):=\oplus_{i=0}^{\infty}S^i (V).
$
Denote the product of elements $v_1,\cdots,v_n\in V$ in $S^n(V)$ by $v_1\odot\cdots\odot v_n$. The degree of $v_1\odot\cdots\odot v_n$ is by definition the sum of the degree of $v_i$. Denote by $\Hom^n(S(V),V)$ the space of degree $n$ linear maps from the graded vector space $S(V)$ to the graded vector space $V$. Obviously, an element $f\in\Hom^n(S(V),V)$ is the sum of $f_i:S^i(V)\lon V$. We will write  $f=\sum_{i=0}^\infty f_i$.
 Set $C^n(V,V):=\Hom^n(S(V),V)$ and
$
C^*(V,V):=\oplus_{n\in\mathbb Z}C^n(V,V).
$

A permutation $\sigma\in\mathbb S_n$ is called an $(i,n-i)$-unshuffle if $\sigma(1)<\cdots <\sigma(i)$ and $\sigma(i+1)<\cdots <\sigma(n)$. If $i=0$ or $n$, we assume $\sigma=\Id$. The set of all $(i,n-i)$-unshuffles will be denoted by $\mathbb S_{(i,n-i)}$. The notion of an $(i_1,\cdots,i_k)$-unshuffle and the set $\mathbb S_{(i_1,\cdots,i_k)}$ are defined analogously.

For a permutation $\sigma\in\mathbb S_n$ and $v_1,\cdots, v_n\in V$,  the Koszul sign $\varepsilon(\sigma;v_1,\cdots,v_n)\in\{-1,1\}$ is defined by
\begin{eqnarray*}
v_1\odot\cdots\odot v_n=\varepsilon(\sigma;v_1,\cdots,v_n)v_{\sigma(1)}\odot\cdots\odot v_{\sigma(n)}.
\end{eqnarray*}
We will abbreviate $\varepsilon(\sigma;v_1,\cdots,v_n)$ to $\varepsilon(\sigma)$.

\begin{defi}{\rm (\cite{LV})}
  Let $(\g=\oplus_{k\in\mathbb Z}\g^k,[\cdot,\cdot],d)$ be a $\dgla$.  A degree $1$ element $\theta\in\g^1$ is called a {\bf Maurer-Cartan element} of $\g$ if it
  satisfies the following {\bf Maurer-Cartan equation}:
  \begin{equation}
  d \theta+\half[\theta,\theta]=0.
  \label{eq:mce}
  \end{equation}
  \end{defi}
A gLa is a dgLa with $d=0$. Then we have

\begin{pro}{\rm (\cite{LV})}
Let $(\g=\oplus_{k\in\mathbb Z}\g^k,[\cdot,\cdot])$ be a $\gla$ and let $\mu\in \g^1$ be a Maurer-Cartan element. Then the map
$$ d_\mu: \g \longrightarrow \g, \ d_\mu(u):=[\mu, u], \quad \forall u\in \g,$$
is a differential on $\g$. For any $v\in \g_1$, the sum $\mu+v$ is a
Maurer-Cartan element of the $\gla$ $(\g,
[\cdot,\cdot])$ if and only if $v$ is a Maurer-Cartan element of the $\dgla$ $(\g, [\cdot,\cdot], d_\mu)$. \label{pp:mce}
\end{pro}

Let $(V;\rho)$ be a representation of a Lie algebra $\g$. Consider the graded vector space
$\huaC^*(V,\g):=\oplus_{k\geq 0}\Hom(\wedge^{k}V,\g).$
Define a skew-symmetric bracket operation
$$\Courant{\cdot,\cdot}: \Hom(\wedge^nV,\g)\times \Hom(\wedge^mV,\g)\longrightarrow \Hom(\wedge^{m+n}V,\g)$$
by taking $f\in\Hom(\wedge^nV,\g)$, $g\in\Hom(\wedge^mV,\g)$ and define
\begin{eqnarray}
&&\nonumber\Courant{f,g}(u_1,u_2,\cdots,u_{m+n})\\
\mlabel{o-bracket}&:=&-\sum_{\sigma\in \mathbb S_{(m,1,n-1)}}(-1)^{\sigma}f(\rho(g(u_{\sigma(1)},\cdots,u_{\sigma(m)}))u_{\sigma(m+1)},u_{\sigma(m+2)},\cdots,u_{\sigma(m+n)})\\
\nonumber&&+(-1)^{mn}\sum_{\sigma\in \mathbb S_{(n,1,m-1)}}(-1)^{\sigma}g(\rho(f(u_{\sigma(1)},\cdots,u_{\sigma(n)}))u_{\sigma(n+1)},u_{\sigma(n+2)},\cdots,u_{\sigma(m+n)})\\
\nonumber&&-(-1)^{mn}\sum_{\sigma\in \mathbb S_{(n,m)}}(-1)^{\sigma}[f(u_{\sigma(1)},\cdots,u_{\sigma(n)}),g(u_{\sigma(n+1)},\cdots,u_{\sigma(m+n)})].
\vspace{-.1cm}
\end{eqnarray}

\begin{pro}\mlabel{pro:gla}{\rm (\cite{TBGS-1})}
  With the above notations, $(\huaC^*(V,\g),\Courant{\cdot,\cdot})$ is a $\gla$. Moreover, its Maurer-Cartan elements are precisely the $\huaO$-operators on $\g$ with respect to the representation $(V;\rho)$.
\end{pro}

Let $T:V\longrightarrow\g$ be an \oop. Since $T$ is a Maurer-Cartan element of the $\gla$ $(\huaC^*(V,\g),\Courant{\cdot,\cdot})$ by Proposition~\ref{pro:gla}, it follows from Proposition~\ref{pp:mce} that
$d_T:=\Courant{T,\cdot}$
 is a graded derivation on the $\gla$ $(\huaC^*(V,\g),\Courant{\cdot,\cdot})$ satisfying $d^2_T=0$.
  Therefore $(\huaC^*(V,\g),\Courant{\cdot,\cdot},d_T)$ is a $\dgla$.
This $\dgla$ controls the deformations of $\huaO$-operators.

\begin{thm}\label{thm:deformation}{\rm (\cite{TBGS-1})}
Let $T:V\longrightarrow\g$ be an $\huaO$-operator on a Lie algebra
$\g$ with respect to a representation $(V;\rho)$. Then for a
linear map $T':V\longrightarrow \g$, $T+T'$ is still
an \oop on the Lie algebra $\g$ with respect to the
representation $(V;\rho)$ if and only if $T'$ is a Maurer-Cartan
element of the $\dgla$
$(\huaC^*(V,\g),\Courant{\cdot,\cdot},d_T)$.
\end{thm}

We next recall the notion of a pre-Lie algebra and the gLa whose Maurer-Cartan elements  characterize pre-Lie algebra structures. We show that there is a close relationship between these two graded Lie algebras.

\begin{defi}\label{de:prelie}{\rm (\cite{Dz})}
  A {\bf pre-Lie algebra} is a pair $(V,\cdot_V)$, where $V$ is a vector space and  $\cdot_V:V\otimes V\longrightarrow V$ is a bilinear multiplication
satisfying that for all $x,y,z\in V$,
$$(x\cdot_V y)\cdot_V z-x\cdot_V(y\cdot_V
z)=(y\cdot_V x)\cdot_V z-y\cdot_V(x\cdot_V z).$$
\end{defi}

Relating an \oop to a pre-Lie algebra, we have
\begin{thm} $($\cite{Bai}$)$
Let $T:V\to \g$ be an $\huaO$-operator on a Lie algebra $\g$ with respect to a representation $(V;\rho)$. Define a multiplication $\cdot_T$ on $V$ by
\begin{equation}
  u\cdot_T v=\rho(Tu)(v),\quad \forall u,v\in V.
\end{equation}
Then $(V,\cdot_T)$ is a pre-Lie algebra, called the {\bf induced
pre-Lie algebra from the $\huaO$-operator $T$}.
\mlabel{thm:opL}
\end{thm}

Let $V$ be a vector space. For $\alpha\in\Hom(\wedge^{n}V\otimes V,V)$ and $\beta\in\Hom(\wedge^{m}V\otimes V,V)$, define  $\alpha\circ\beta\in\Hom(\wedge^{n+m}V\otimes V,V)$ by
\begin{eqnarray}
\nonumber&&(\alpha\circ\beta)(u_1,\cdots,u_{m+n+1})\\
\mlabel{eq:pLbrac}&:=&\sum_{\sigma\in\mathbb S_{(m,1,n-1)}}(-1)^{\sigma}\alpha(\beta(u_{\sigma(1)},\cdots,u_{\sigma(m+1)}),u_{\sigma(m+2)},\cdots,u_{\sigma(m+n)},u_{m+n+1})\\
&&+(-1)^{mn}\sum_{\sigma\in\mathbb S_{(n,m)}}(-1)^{\sigma}\alpha(u_{\sigma(1)},\cdots,u_{\sigma(n)},\beta(u_{\sigma(n+1)},\cdots,u_{\sigma(m+n)},u_{m+n+1})).
\nonumber
\end{eqnarray}
Then the graded vector space
$\CV^*(V,V):=\oplus_{k\geq 0} \Hom(\wedge^kV\ot V,V)$ is a gLa when it is equipped with the {\bf Matsushima-Nijenhuis bracket} $[\cdot,\cdot]^C$~ \cite{CL,Nijenhuis,WBLS}:
\begin{eqnarray}
[\alpha,\beta]^C:=\alpha\circ\beta-(-1)^{mn}\beta\circ\alpha,
\quad \forall \alpha\in\Hom(\wedge^{n}V\otimes
V,V),\beta\in\Hom(\wedge^{m}V\otimes V,V).
\mlabel{eq:glapL}
\end{eqnarray}

\begin{rmk}
In fact, $\alpha\in\Hom(V\otimes V,V)$ defines a pre-Lie algebra
structure on $V$ if and only if $[\alpha,\alpha]^C=0,$ that is,
$\alpha$ is a Maurer-Cartan element  of the $\gla$
$(\CV^*(V,V),[\cdot,\cdot]^C)$. \label{rk:plmc}
\end{rmk}

Define a linear map $\Phi:\Hom(\wedge^kV,\g)\lon\Hom(\wedge^{k}V\otimes V,V), k\geq 0$ by
\begin{eqnarray}\label{eq:phi}
\Phi(f)(u_1,\cdots,u_k,u_{k+1})=\rho(f(u_1,\cdots,u_k))(u_{k+1}),\,\,\,\,\forall
f\in\Hom(\wedge^kV,\g), u_1,\cdots,u_{k+1}\in V.
\mlabel{eq:defiphi}
\end{eqnarray}

\begin{thm}{\rm (\cite{TBGS-1})}
Let $(V;\rho)$ be a representation of a Lie algebra $\g$. Then $\Phi$ is a homomorphism  of gLas from $(\huaC^*(V,\g),\Courant{\cdot,\cdot})$ to
$(\CV^*(V,V),[\cdot,\cdot]^C)$.
\end{thm}

\begin{rmk}\label{rmk:O-pre}
As a direct consequence of the above proposition, the
Maurer-Cartan elements in the first graded Lie algebra are sent to
those in the second graded Lie algebra. Thus by
Proposition~\ref{pro:gla} and Remark~\ref{rk:plmc}, the \oops
are sent to pre-Lie algebra structures on $V$. Furthermore, two \oops  are sent to the same pre-Lie algebra if and only if they are in the same fiber of $\Phi$. This lifts the connection from $\calo$-operators to pre-Lie algebras in Theorem~\ref{thm:opL} to the level of deformations.
\end{rmk}

\section{Homotopy $\huaO$-operators on symmetric  graded Lie algebras}
\label{sec:hom}

We first recall the symmetric generalizations the notions of gLas and dgLas~\cite{ALN}.

\begin{defi}\label{symmetric-lie}
A {\bf symmetric  graded Lie algebra (sgLa)} is a $\mathbb Z$-graded vector space $\g=\oplus_{k\in\mathbb Z}\g_k$ equipped with a bilinear bracket $[\cdot,\cdot]_\g:\g\otimes \g\lon \g$  of degree $1$ such that
\begin{itemize}
\item[\rm(i)] {\bf (graded symmetry)} $[x,y]_\g=(-1)^{xy}[y,x]_\g,$
\item[\rm(ii)] {\bf (graded Leibniz rule)} $[x,[y,z]_\g]_\g=(-1)^{x+1}[[x,y]_\g,z]_\g+(-1)^{(x+1)(y+1)}[y,[x,z]_\g]_\g$.
\end{itemize}
Here $x,y,z$ are homogeneous elements in $\g$, which also denote their degrees when in exponent.
\end{defi}

\begin{defi}\label{symmetric-dg-lie}
A {\bf symmetric differential graded Lie algebra (sdgLa)} is a symmetric  graded Lie algebra $(\g,[\cdot,\cdot]_\g)$ equipped with a linear map $d:\g\lon\g$ of degree $1$ such that
\begin{eqnarray}
d[x,y]_\g=-[dx,y]_\g-(-1)^x[x,dy]_\g.
\end{eqnarray}
\end{defi}

We also recall the notion of the suspension and desuspension operators. Let $V=\oplus_{i\in\mathbb Z} V^i$ be a graded vector space, we define the {\bf suspension operator} $s:V\mapsto sV$ by assigning  $V$ to the graded vector space $sV=\oplus_{i\in\mathbb Z}(sV)^i$ with $(sV)^i:=V^{i-1}$.
There is a natural degree $1$ map $s:V\lon sV$ that is the identity map of the underlying vector space, sending $v\in V$ to its suspended copy $sv\in sV$. Likewise, the {\bf desuspension operator} $s^{-1}$ changes the grading of $V$ according to the rule $(s^{-1}V)^i:=V^{i+1}$. The  degree $-1$ map $s^{-1}:V\lon s^{-1}V$ is defined in the obvious way.

Let $(\g,[\cdot,\cdot]_\g)$ and $(\g',[\cdot,\cdot]_{\g'})$ be \sgla's. A {\bf homomorphism} from $\g$ to $\g'$ is a linear map $\phi:\g\lon \g'$ of degree $0$ such that
\begin{eqnarray*}
\phi([x_1,x_2]_\g)=[\phi(x_1),\phi(x_2)]_{\g'},\quad\forall x_1,x_2\in \g.
\end{eqnarray*}

\begin{defi}
A {\bf representation} of an \sgla~$(\g,[\cdot,\cdot]_\g)$ on a graded vector space $V$ is a homomorphism of graded vector spaces $\rho:\g\lon \gl(V)$ of degree $1$ such that $s^{-1}\circ\rho:\g\lon s^{-1}\gl(V)$ is an \sgla~ homomorphism.
\end{defi}

Now we are ready to give the main notion of this section.

\begin{defi}{\rm (\cite{TBGS-2})}
Let  $\rho$ be a representation of an \sgla~ $(\g,[\cdot,\cdot]_\g)$ on a graded vector space $V$. A degree $0$ element $T=\sum_{i=0}^{+\infty}T_i\in \Hom(S(\h),\g)$ with $T_i\in \Hom(S^i(\h),\g)$ is called a  {\bf homotopy $\huaO$-operator} on an \sgla~ $(\g,[\cdot,\cdot]_\g)$ with respect to the representation $\rho$ if the following generalized Rota-Baxter identities hold for all $p\geq 0$ and all homogeneous elements $v_1,\cdots,v_p\in V$,
\begin{eqnarray}
\label{homotopy-rota-baxter-o}\nonumber&&\sum_{k+l=p+1}\sum_{\sigma\in \mathbb S_{(l,1,p-l-1)}}\varepsilon(\sigma)T_{k-1}\Big(\rho\big(T_l(v_{\sigma(1)},\cdots,v_{\sigma(l)})\big)v_{\sigma(l+1)},v_{\sigma(l+2)},\cdots,v_{\sigma(p)}\Big)\\
&=&\frac{1}{2}\sum_{k+l=p+1}\sum_{\sigma\in \mathbb S_{(k-1,l)}}\varepsilon(\sigma)[T_{k-1}(v_{\sigma(1)},\cdots,v_{\sigma(k-1)}),T_l(v_{\sigma(k)},\cdots,v_{\sigma(p)})]_\g.
\end{eqnarray}
\label{de:homoop}
\end{defi}

\begin{rmk}
The linear map $T_0$ is just an element $\Omega\in \g^0$. Below are the generalized Rota-Baxter identities for $p=0,1,2:$
\begin{eqnarray*}
\label{weak-1}[\Omega,\Omega]_\g&=&0,\\
\label{weak-2}T_1(\rho(\Omega)v_1)&=&[\Omega,T_1(v_1)]_\g,\\
\label{weak-3}[T_1(v_1),T_1(v_2)]_\g&=&T_1\Big(\rho(T_1(v_1))v_2+(-1)^{v_1v_2}\rho(T_1(v_2))v_1\Big)\\
&&+T_2(\rho(\Omega)v_1,v_2)+(-1)^{v_1v_2}T_2(\rho(\Omega)v_2,v_1)-[\Omega,T_2(v_1,v_2)]_\g.
\end{eqnarray*}
\end{rmk}

\begin{rmk}
  If the \sgla~  reduces to a Lie algebra and the action reduces to a representation of a Lie algebra on a vector space, the above definition reduces to the definition of an $\huaO$-operator on a Lie algebra.
\end{rmk}

\begin{rmk}
A {\bf homotopy Rota-Baxter operator} $R=\sum_{i=0}^{+\infty}R_i\in \Hom(S(\g),\g)$   on an \sgla~ $(\g,[\cdot,\cdot]_\g)$ is a homotopy $\huaO$-operator with respect to the adjoint representation $\ad$.
If moreover the \sgla~  reduces to a Lie algebra, then the resulting linear operator $R:\g\longrightarrow \g$ is a { Rota-Baxter operator}.
\end{rmk}

In the sequel, we construct a \gla~ and show that  homotopy $\huaO$-operators can be characterized as its Maurer-Cartan elements to justify our definition of homotopy $\huaO$-operators.

Let $\rho$ be a representation of an \sgla~ $(\g,[\cdot,\cdot]_\g)$ on a graded vector space $V$. Consider the graded vector space
$C^*(V,\g):=\oplus_{n\in\mathbb Z}\Hom^n(S(V),\g)$.
Define a graded bracket operation
 $$\Courant{\cdot,\cdot}: \Hom^m(S(V),\g)\times \Hom^n(S(V),\g)\longrightarrow \Hom^{m+n+1}(S(V),\g)$$
by
\begin{eqnarray}
\nonumber&&\Courant{f,g}_p( v_1,\cdots,v_{p})\\
\label{graded-Lie}&=&-\sum_{k+l=p+1}\sum_{\sigma\in \mathbb S_{(l,1,p-l-1)}}\varepsilon(\sigma)f_{k-1}\Big(\rho\big(g_l(v_{\sigma(1)},\cdots,v_{\sigma(l)})\big)v_{\sigma(l+1)},v_{\sigma(l+2)},\cdots,v_{\sigma(p)}\Big)\\
\nonumber&&+(-1)^{(m+1)(n+1)}\sum_{k+l=s+1}\sum_{\sigma\in \mathbb S_{(k-1,1,p-k)}}\varepsilon(\sigma)g_l\Big(\rho\big(f_{k-1}(v_{\sigma(1)},\cdots,v_{\sigma(k-1)})\big)v_{\sigma(k)}, v_{\sigma(k+1)},\cdots,v_{\sigma(p)}\Big)\\
\nonumber&&-\sum_{k+l=p+1}\sum_{\sigma\in \mathbb S_{(k-1,l)}}(-1)^{n(v_{\sigma(1)}+\cdots+v_{\sigma(k-1)})+m+1}\varepsilon(\sigma)[f_{k-1}(v_{\sigma(1)},\cdots,v_{\sigma(k-1)}),g_l(v_{\sigma(k)},\cdots,v_{\sigma(p)})]_\g
\end{eqnarray}
for all  $f=\sum_if_i\in\Hom^m(S(V),\g)$, $g=\sum g_i\in \Hom^n(S(V),\g)$  with  $f_i, g_i\in\Hom(S^i(V),\g)$
 and $v_1,\cdots, v_{p} \in V.$ Here we write $\Courant{f,g}=\sum_i\Courant{f,g}_i$ with $\Courant{f,g}_i\in\Hom(S^i(V),\g)$.

\begin{thm}\label{dgla-deforamtion-homotopy}{\rm (\cite{TBGS-2})}
Let $\rho$ be a representation of an \sgla~ $(\g,[\cdot,\cdot]_\g)$ on a graded vector space $V$. Then $(sC^*(V,\g),\Courant{\cdot,\cdot})$ is a \gla.
\end{thm}

Homotopy  $\huaO$-operators can be characterized as Maurer-Cartan elements of the above \gla. Note that an element $T=\sum_{i=0}^{+\infty}T_i\in \Hom(S(V),\g)$ is of degree $0$ if and only if the corresponding element $T\in s\Hom(S(V),\g)$ is of degree $1$.

\begin{thm}\label{hmotopy-o-operator-dgla}{\rm (\cite{TBGS-2})}
Let $\rho$ be a representation of an \sgla~ $(\g,[\cdot,\cdot]_\g)$ on a graded vector space $V$. A degree $0$ element $T=\sum_{i=0}^{+\infty}T_i\in \Hom(S(V),\g)$ is a homotopy $\huaO$-operator on $\g$ with respect to the representation $\rho$ if and only if $T=\sum_{i=0}^{+\infty}T_i$ is a Maurer-Cartan element of the \gla $(sC^*(V,\g),\Courant{\cdot,\cdot})$.
\end{thm}
The notion of a pre-Lie$_\infty$-algebra was introduced in \cite{CL} as the homotopy of the pre-Lie algebra. See \cite{Mer} for applications of pre-Lie$_\infty$-algebras in geometry.
\begin{defi}
  A {\bf pre-Lie$_\infty$-algebra} is a graded vector space $V$ equipped with a collection of linear maps $\oprn_k:\otimes^k V\lon V, k\ge 1,$  of degree $1$ with the property that, for any homogeneous elements $v_1,\cdots,v_n\in V$, we have
\begin{itemize}\item[\rm(i)]
{\bf (graded symmetry)} for every $\sigma\in\mathbb S_{n-1}$,
\begin{eqnarray*}
\oprn_n(v_{\sigma(1)},\cdots,v_{\sigma(n-1)},v_n)=\varepsilon(\sigma)\oprn_n(v_1,\cdots,v_{n-1},v_n),
\end{eqnarray*}
\item[\rm(ii)] for all $n\ge 1$,
\begin{eqnarray*}
&&\sum_{i+j=n+1\atop i\ge1,j\geq2}\sum_{\sigma\in\mathbb S_{(i-1,1,j-2)}}\varepsilon(\sigma)\oprn_j(\oprn_i(v_{\sigma(1)},\cdots,v_{\sigma(i-1)},v_{\sigma(i)}), v_{\sigma(i+1)},\cdots,v_{\sigma(n-1)},v_{n})\\
&&+\sum_{i+j=n+1\atop i\geq1,j\geq1}\sum_{\sigma\in\mathbb S_{(j-1,i-1)}}(-1)^{\alpha}\varepsilon(\sigma)\oprn_j(v_{\sigma(1)},\cdots,v_{\sigma(j-1)},\oprn_i( v_{\sigma(j)},\cdots,v_{\sigma(n-1)},v_{n}))=0,
\end{eqnarray*}
where $\alpha=v_{\sigma(1)}+v_{\sigma(2)}+\cdots+v_{\sigma(j-1)}$.
\end{itemize}
\end{defi}

Pre-Lie$_\infty$-algebras can be characterized as the Maurer-Cartan elements of a dgLa.

\begin{pro} {\rm (\cite{CL})}
Let $V$ be a $\mathbb Z$-graded vector space. Denote by
$$
\bar{\huaC}^n(V,V):=s\Hom^{n-1}(S(V),s^{-1}\gl(V))=\Hom^n(S(V),\gl(V)),\quad \bar{\huaC}^*(V,V):=\oplus_{n\in\mathbb Z} \bar{\huaC}^n(V,V).
$$
Then $(\bar{\huaC}^*(V,V),[\cdot,\cdot]^c)$ is a gla  equipped with the graded Lie bracket
 $$[\cdot,\cdot]^c:\bar{\huaC}^m(V,V)\times \bar{\huaC}^n(V,V)\longrightarrow \bar{\huaC}^{m+n}(V,V)$$
whose details is referred to the original literature.

Moreover, $L=\sum_{i=0}^{\infty}L_i\in\Hom^1(S(V),\gl(V))$ defines a pre-Lie$_\infty$-algebra structure by
\begin{eqnarray}
\label{homotopy-pre-lie}\oprn_k(v_1,\cdots,v_k):=L_{k-1}(v_1,\cdots,v_{k-1})v_k,\quad\forall v_1\cdots,v_k\in V,
\end{eqnarray}
on the graded vector space $V$ if and only if  $L=\sum_{i=0}^{\infty}L_i$ is a  Maurer-Cartan element of the \gla~ $(\bar{\huaC}^*(V,V),[\cdot,\cdot]^c)$.
\end{pro}

In the above result, if the graded vector space $V$ reduces to a usual vector space, we obtain the {\bf Matsushima-Nijenhuis bracket} $[\cdot,\cdot]^C$ given by \eqref{eq:glapL}. See \cite{CL,Nijenhuis,WBLS} for more details.

Finally, as promised, we lift the connection from $\calo$-operator to pre-Lie algebras (Theorem~\ref{thm:opL}) to the level of homotopy. We first establish the relationship between the two gLas that give Maurer-Cartan characterizations of homotopy $\huaO$-operators and pre-Lie$_\infty$-algebras respectively. Define a graded linear map $\Psi:sC^*(V,\g)\lon\bar{\huaC}^*(V,V)$ of degree $0$ by
$$
\Psi(f)=s^{-1}\circ\rho\circ f,\quad \forall f\in\Hom^m(S(V),\g).
$$
Therefore, we have $\Psi(f)_k=s^{-1}\circ\rho\circ f_k.$

\begin{thm}\label{homo-dg-lie}{\rm (\cite{TBGS-2})}
Let $\rho$ be a representation of an \sgla~ $(\g,[\cdot,\cdot]_\g)$ on a graded vector space $V$.   Then $\Psi$ is a
homomorphism of gLas from $(sC^*(V,\g),\Courant{\cdot,\cdot})$  to
$(\bar{\huaC}^*(V,V),[\cdot,\cdot]^c)$.
\end{thm}

By this theorem, we can obtain a pre-Lie$_\infty$-algebra from a homotopy $\huaO$-operator.

\begin{cor}\label{homotopy-RB-homotopy-post-lie}
Let $T=\sum_{i=0}^{+\infty}T_i\in \Hom^0(S(V),\g)$ be a homotopy $\huaO$-operator on an \sgla~ $(\g,[\cdot,\cdot]_\g)$ with respect to a representation $\rho:\g\lon\gl(V)$. Then $(V,\{\oprn_k\}_{k=1}^\infty)$ is a pre-Lie$_\infty$-algebra, where  $\oprn_k:\otimes^k V\lon V$ $(k\ge 1)$  are   linear maps  of degree $1$ defined by
\begin{eqnarray}
\label{homotopy-o-to-homotopy-post-lie}\oprn_k(v_1,\cdots,v_k):=\rho(T_{k-1}(v_1,\cdots,v_{k-1})\big)v_k,\quad \forall v_1\cdots,v_k\in V.
\end{eqnarray}
\end{cor}

\noindent{\bf Acknowledgements. } This research is supported by
NSFC (11922110,  11425104, 11771190, 11931009). C. Bai
is also supported by the Fundamental Research Funds for the
Central Universities and Nankai ZhiDe Foundation.

\end{document}